\documentclass[a4paper]{article}
\usepackage{amssymb,amsmath}

\newcommand{\rcdbpd}[2]{\frac{\partial {#1} }{\partial {#2} } }
\newcommand{\rcdbdpd}[2]{\dfrac{\partial {#1} }{\partial {#2} } }
\newtheorem{example}{Example}
\newtheorem{theorem}{Theorem}

\begin{document}

\begin{center}

{\bf\Large Jet Bundles and the Formal Theory of \\
Partial Differential Equations}

\vspace{0.4cm}

Richard Baker\footnote{e-mail: \texttt{rcdb100@mrao.cam.ac.uk}}  and
Chris Doran\footnote{e-mail: \texttt{c.doran@mrao.cam.ac.uk},  
\texttt{http://www.mrao.cam.ac.uk/$\sim$cjld1/}} 

\vspace{0.4cm}

Astrophysics Group, Cavendish Laboratory, Madingley Road, \\
Cambridge CB3 0HE, UK.

\vspace{0.4cm}

\begin{abstract}
Systems of partial differential equations lie at the heart of physics.
Despite this, the general theory of these systems has remained rather
obscure in comparison to numerical approaches such as finite element
models and various other discretisation schemes. There are, however,
several theoretical approaches to systems of PDEs, including schemes
based on differential algebra and geometric approaches including the
theory of exterior differential systems~\cite{amp} and the so-called
``formal theory''~\cite{th:seiler} built on the jet bundle formalism.
This paper is a brief introduction to jet bundles focusing on the
completion of systems to equivalent involutive systems for which power
series solutions may be constructed order by order. We will not consider
the mathematical underpinnings of involution (which lie in the theory of
combinatorial decompositions of polynomial modules~\cite{calmet,seiler})
nor other applications of the theory of jet bundles such as the theory
of symmetries of systems of PDEs~\cite{olver} or discretisation schemes
based on discrete approximations to jet bundles~\cite{marsden}.
\end{abstract}

\end{center}

\section{Fibre Bundles and Sections}

A bundle is a triple $\left( M, X, \pi \right)$, where $M$ is a
manifold called the total space, $X$ is a manifold called the base
space, and $\pi : M \rightarrow X$ is a continuous surjective mapping
called the projection. Where no confusion can arise, we shall often
find it convenient to denote the bundle either by its total space or
its projection.  A trivial bundle is a bundle whose total space $M$ is
homeomorphic to $X\times U$, where $U$ is a manifold called the
fibre. A bundle which is locally a trivial bundle is called a fibre
bundle (an example of a fibre bundle which is not trivial is the
M\"obius band). In the following we shall only be concerned with
trivial bundles.

We will denote the coordinates on the base manifold $x = \{x^i, i=1,
\dots ,p\}$ where $p$ is the dimension of the base manifold and the
coordinates on the fibre $u = \{u^\alpha, \alpha=1, \dots, q\}$ where
$q$ is the dimension of the fibre (when we consider jet bundles, we will
need to extend this notation slightly). In the case of a trivial bundle,
or in a local coordinate patch on a fibre bundle, the projection takes
the simple form
\[
\pi : \left\{
\begin{array}{l}
X \times U \rightarrow X \\
(x,u) \mapsto (x)
\end{array}
\right.
\]
A section of a fibre bundle is a map
\[
\Phi_{f} : \left\{ \begin{array}{ll}
                X \rightarrow X \times U  \\
                x \mapsto \left( x, f \left( x \right) \right)
              \end{array}
       \right.
\]
such that $\pi \circ \Phi_{f}$ is the identity map on $X$. In other
words, a section assigns to each point in $X$ a point in the fibre over
that point. The graph of the function $f(x)$ is:
\[
\Gamma_{f} = \left\{ \left( x,f(x) \right) : x \in \Omega
\right\} \subset X \times U
\]
where $\Omega$ is the domain of definition of $f$. We will find it
convenient to refer to sections, functions and graphs interchangably.

\section{Jet Bundles}

We define a multi-index $J$ as a $p$-tuple $[j_1,j_2,\dots,j_p]$ with
$j_i \in \mathbb{N}$. The order of the multi-index $J$, denoted $|J|$ is
given by the sum of the $j_i$. We will often find it more convenient to
use a repeated-index notation for $J$. In this notation $J$ is
represented by a string of $|J|$ independent coordinate labels, with
$j_i$ copies of the $i$-th coordinate label. For example, if $p=3$ and
the coordinates are labelled then $x$, $y$ and $z$ then the second order
multi-indices in repeated index notation are $xx$, $xy$, $xz$, $yy$,
$yz$ and $zz$. We introduce the special notation $J,i$ where $i$ is an
independent coordinate label for the multi-index given by
$[j_1,\dots,j_i+1,\dots,j_p]$. For example, $xyy,x = xxyy$.

If our independent variables are $x^i$ and our dependent variables are
$u^\alpha$ then we introduce jet variables $u^\alpha_J$ where $J$ is a
multi-index. Notice that we can put the jet variables of order $n$ in
one-to-one correspondance with the derivatives of the dependent
variables of order $n$. We will later introduce further structures that
enable us to make a full correspondance between jet variables and
derivatives.

Associated with these jet variables we introduce a set of Euclidean
spaces $U_i$, whose coordinates are $u^\alpha_J$ with $|J|=i$. We call
the space $M^{(1)} = X \times U \times U_1$ the first order jet bundle
over the space $M = X \times U$.  We now introduce the notation
\[
U^{(n)} = U \times U_1 \times \dots \times U_n
\]
and call the space $M^{(n)} = X \times U^{(n)}$ the $n$-th order jet
bundle over $M$.

\begin{example}
Let $p=2$ and $q=1$. Label the independent variables $x$ and $y$ and 
the
dependent variable $u$. The first order jet bundle, $M^{(1)}$, then has
coordinates $(x,y,u,u_x,u_y)$, the second order jet bundle, $M^{(2)}$,
has coordinates $(x,y,u,u_x,u_y,u_{xx},u_{xy},u_{yy})$ and so on.
\end{example}

We will often consider a jet bundle as a bundle over a lower order jet
bundle. We denote the natural projection between the $(m+n)$-th order
jet bundle and the $n$-th order jet bundle as
\[
  \pi^{m+n}_{n} : M^{(m+n)} \rightarrow M^{(n)}
\]
Note that although $M^{(m+n)}$ is a bundle over $M^{(n)}$ it is not a
jet bundle over $M^{(n)}$, but rather a subset of such a bundle.

\section{Differential Functions and Formal Derivatives}

A differential function is a smooth, real-valued function defined on
$M^{(n)}$. We denote the algebra of differential functions defined on
$M^{(n)}$ by $\mathcal{A}_{(n)}$. If $F \in \mathcal{A}_{(n)}$ then $F
\in \mathcal{A}_{(m+n)}$ too, as the coordinates on $M^{(n)}$ are a
subset of the coordinates on $M^{(m+n)}$. If the lowest order space on
which $F$ is defined is $M^{(n)}$ then we will say that $F$ is an
$n$-th order differential function.  These will be used to describe
sections of $M$ and differential equations.

The most fundamental maps between lower order and higher order jet
bun\-dles are provided by formal derivatives. A formal derivative operator
$D_i$, called the formal derivative with respect to $x^i$, maps each
differential function $F \in \mathcal{A}_{(n)}$ to a differential
function $D_iF \in \mathcal{A}_{(n+1)}$ via
\[
D_iF = \rcdbpd{F}{x^i} + \sum_{\alpha=1}^q \sum_J \rcdbpd{F}{u^\alpha_J}
\, u^\alpha_{J,i}
\]

It is convenient to extend the notation for formal and
partial derivatives to encompass our multi-index notation:
\[
D_J = (D_1)^{j_1} (D_2)^{j_2} \dots (D_p)^{j_p}
\]
and similarly for partial derivatives. Clearly if $F \in
\mathcal{A}_{(n)}$ then $D_JF \in \mathcal{A}_{(n+|J|)}$.

\section{Prolongation of Sections}

If $\Gamma_f \subset M$ is a section defined by $u = f(x)$ where $f$ is
a smooth function of $x$ then we can use the formal derivative to
prolong it to a section $\Gamma_f^{(1)} \subset M^{(1)}$. The equations
defining $\Gamma_f^{(1)}$ are simply $u^\alpha = f^a(x)$ and the
equations found by applying each of the $D_i$ to $u^\alpha =
f^\alpha(x)$:
\[
\Gamma_f^{(1)} = \left\{
  \begin{array}{ccc}
    u^\alpha & = & f^\alpha(x) \\
    u^\alpha_i & = & \partial_i f^\alpha(x)
  \end{array}
\right.
\]
Similarly we can use the $D_i$ multiple times to prolong $\Gamma_f$ to a
section $\Gamma_f^{(n)} \subset M^{(n)}$ defined by the equations
\[
u_J^\alpha = \partial_J f^\alpha(x)
\]
where $J$ ranges over all multi-indices such that $0 \leqslant |J|
\leqslant n$. We will sometimes talk about the $n$-th prolongation of a
function $f(x)$, and write this prolongation as $f^{(n)}(x)$.

\section{Differential Equations and Solutions}

We intend to view systems of PDEs as geometric objects. In keeping
with this programme, we will simply call such a system a ``differential
equation''. An $n$-th order differential equation, $I_\Delta$, is a
fibred submanifold of $M^{(n)}$. The differential equation is often
stated as the kernel of a set of differential functions $\Delta_\nu \in
\mathcal{A}_{(n)}$:
\[
\begin{array}{lr}
\Delta_\nu(x,u^{(n)}) = 0, & \nu = 1,\dots,l
\end{array}
\]
We can map any system of partial differential equations onto such a
submanifold, simply by replacing all of the derivatives of dependent
variables by the corresponding jet variables. Indeed we have chosen our
notation in such a way that this process is entirely transparent.

\begin{example}
The two dimensional wave equation
\[
\frac{\partial^2 u}{\partial t^2} - \frac{\partial^2 u}{\partial x^2} -
\frac{\partial^2 u}{\partial y^2} = 0
\]
maps onto the submanifold of $M^{(2)}$ (with the obvious coordinates)
determined by
\[
u_{tt} - u_{xx} - u_{yy} = 0
\]
\end{example}

In keeping with our compact notations $x$, $u$ and $f$, we shall write
\[
\Delta\left(x,u^{(n)}\right) =
\left(\Delta_1(x,u^{(n)}),\dots,\Delta_l(x,u^{(n)})\right) \in
\mathcal{A}_{(n)}^l
\]
$\Delta$ is therefore a map from $M^{(n)}$ to $\mathbb{R}^l$. The
differential equation $I_\Delta$ is the submanifold in which the map
$\Delta$ vanishes:
\[
I_\Delta = \left\{ (x,u^{(n)}) : \Delta(x,u^{(n)}) = 0
\right\}
\]
A smooth solution of $I_\Delta \subset M^{(n)}$ is a smooth function
$f(x)$ such that
\[
\Delta(x,f^{(n)}(x)) = 0
\]
or, in terms of our geometric formulation
\[
\Gamma_f^{(n)} \subset I_\Delta
\]
Note that not every section of $M^{(n)}$ which lies entirely within
$I_\Delta$ is a prolongation of a section of $M$ - the prolongation of a
section of $M$ automatically respects the correspondence between jet
variables and derivatives, whereas an arbitrary section of $M^{(n)}$
does not.

\section{Prolongation and Projection of Differential \\ Equations}

The $k$-th prolongation of the differential equation
\[
I_\Delta = \left\{ (x,u^{(n)}) : \Delta_\nu(x,u^{(n)})=0 \right\}
\subset M^{(n)}
\]
is
\[
I^{(k)}_\Delta = \left\{ (x,u^{(n+k)}) :
(D_J\Delta_\nu)(x,u^{(n+k)})=0 \right\} \subset M^{(n+k)}
\]
where $J$ runs over all multi-indices up to order $k$.

Differential equations may be projected along the fibers onto lower
order jet bundles. In general, this is a complicated procedure in local
coordinates, but it is much easier if the differential equation is known
to be the prolongation of a lower order system. For the remainder of
this paper this will always be the case.

\begin{example}
Let the differential equation $\mathcal{I} \subset M^{(2)}$ be defined by
\[
\mathcal{I} : \left\{
  \begin{array}{r}
    u_{zz} + u_{xy} + u = 0 \\
    u_x - u = 0 \\
    u_y - u^2 = 0 \\
  \end{array}
\right.
\]
The projection of $\mathcal{I}$ into $M^{(1)}$ is
\[
\pi^2_1 (\mathcal{I}) : \left\{
  \begin{array}{r}
    u_x - u = 0 \\
    u_y - u^2 = 0 \\
  \end{array}
\right.
\]

\end{example}

\section{Power Series Solutions}

A smooth solution of $\mathcal{I} \subset M^{(n)}$ in a neighbourhood of
$x_0$ may be written as the power series
\[
u^\alpha = f^\alpha (x) = \sum_{|J|=0}^{\infty} \frac{a^\alpha_J}{J!} 
\left(
x-x_0 \right)^J
\]
for some constants $a^\alpha_J$. Here $J! = j_1! \; j_2! \; \dots \;
j_p!$ and
\[
  \left( x-x_0 \right)^J = \prod_{i=1}^p \left( x^i-x^i_0 \right)^{j_i}
\]
All we have to do is choose the values of the $a^\alpha_J$ so that
$\Gamma_f^{(n)} \subset \mathcal{I}$. It will, however, prove to be
easiest to fix the $a^\alpha_J$ by using the condition that
$\Gamma_f^{(n+k)} \subset \mathcal{I}^{(k)}$ for all $k$ and working
entirely at the point $x_0$. By applying the formal derivative
repeatedly to the power series we can obtain power series expressions
for each of the $u^\alpha_J(x)$. Evaluating each of these at $x_0$ shows
that $a^\alpha_J = u^\alpha_J|_{x_0}$. Therefore we require that
\[
D_J \Delta_\nu(x_0,a^\alpha_J) = 0
\]
for all $J$. We have exchanged the solution of a set of partial
differential equations for the solution of an infinite number of
algebraic equations.

For a class of systems known as formally integrable differential
equations we can construct a power series solution order by order. We
first substitute the general form for the power series into each of the
equations and evaluate at $x=x_0$. This gives us a set of algebraic
equations for the $a^\alpha_J(x)$ with $|J| \leqslant n$. We then make a
partition of the jet variables into parametric derivatives whose values
we can choose and principal derivatives whose values are then fixed by
the system, and solve for the latter in terms of the former. We then
prolong the differential equation and repeat the process. This time the
equations of order less than $n+1$ will automatically be satisfied by
the previously chosen constants, and we will be left with a new set of
equations $D_i \Delta_\nu (x_0,a^\alpha_J) = 0$ for $|J|=n+1$. The
nature of the formal derivative means that these equations will be
linear. We may repeat the procedure to calculate ever higher terms in
the power series.

\section{Integrability Conditions}

If a differential equation is not formally integrable the solutions are
subject to constraints which we call integrability conditions. These are
extra equations that are differential rather than algebraic consequences
of the equations $\Delta = 0$. In other words, projecting the
prolongation of a differential equation may not return the original
equation but only a proper subset thereof:
\[
\pi^{n+k}_n ( \mathcal{I}^{(k)}_\Delta ) \varsubsetneq
\mathcal{I}_\Delta \subset M^{(n)}
\]
and so the order by order construction of a power series solution will
be disrupted.  To streamline the notation, we will write the $j$-th
projection of the $k$-th prolongation of $\mathcal{I}$ as
$\mathcal{I}^{(k)}_j$. For example, the expression above may be
rewritten $\mathcal{I}^{(k)}_k
\varsubsetneq \mathcal{I}$.

Integrability conditions arise in two ways: through the differentiation
of equations with order less than $n$ in $\mathcal{I} \subset M^{(n)}$,
and through the effects of cross-derivatives, as shown in the following
example:

\begin{example}
\label{ex:integrability}
Let $p=3$ with coordinates $x$, $y$ and $z$, and $q=1$ with coordinate
$u$, and consider the differential equation
\[
\mathcal{I} : \left\{
  \begin{array}{r}
    u_z + y u_x = 0 \\
    u_y = 0
  \end{array}
\right.
\]
which prolongs to
\[
\mathcal{I}^{(1)} : \left\{
  \begin{array}{r}  
    u_{xz} + y u_{xx} = 0 \\
    u_{yz} + u_x + y u_{xy} = 0 \\
    u_{zz} + y u_{xz} = 0 \\   
    u_{xy} = 0 \\
    u_{yy} = 0 \\
    u_{yz} = 0 \\
    u_z + y u_x = 0 \\
    u_y = 0
  \end{array}
\right.
\]
We see that the equations $u_{yz} = 0$ and $u_{xy} = 0$ substituted into
the equation $u_{yz} + u_x + y u_{xy} = 0$ imply that $u_x = 0$. This is
a first order equation and so forms part of $\mathcal{I}^{(1)}_1$ on
projection. Hence $\mathcal{I}^{(1)}_1 \varsubsetneq \mathcal{I}$.

\end{example}

A differential equation $\mathcal{I} \subset M^{(n)}$ is formally
integrable if $\mathcal{I}^{(k+1)}_1 = \mathcal{I}^{(k)}$ for all $k$.
Notice that to check for formal integrability requires an infinite
number of operations, for integrability conditions may in general arise
after an arbitrarily large number of prolongations.

\section{Involutive Differential Equations}

We now turn to the consideration of a subset of formally integrable
differential equations known as involutive equations. Two facts make
this class of equations interesting and useful. Firstly, it is possible
to determine whether a given differential equation is involutive using
only a finite number of operations. Secondly, for any differential
equation it is possible to produce an involutive equation with the same
solution space using only a finite number of operations.

There is a more systematic method for determining the integrability
conditions that arise upon a single prolongation and projection. Let us
look at the Jacobi matrix of $\mathcal{I}^{(1)} \subset M^{(n+1)}$. This
matrix can be divided into four blocks:
\[
\left( \,
\begin{array}{|c|c|}
\hline
  \begin{array}{lr}
    \rcdbdpd{D_i\Delta_\nu}{u^\alpha_J}, & 0 \leqslant |J| \leqslant n
  \end{array}
&
  \begin{array}{lr}
    \rcdbdpd{D_i\Delta_\nu}{u^\alpha_J}, & |J| = n+1
  \end{array}
\\
\hline
  \begin{array}{lr}
    \rcdbdpd{\Delta_\nu}{u^\alpha_J}, & 0 \leqslant |J| \leqslant n
  \end{array}
&
0 \\
\hline
\end{array}
\, \right)
\]
We order the columns according to increasing $|J|$, and within each
order by the first non-vanishing component of J (we will call this
component the class of J).

When we project the $\mathcal{I}^{(1)}$ into $M^{(n)}$ to form
$\mathcal{I}^{(1)}$ we must include only those equations that are
independent of the $u^\alpha_J$ with $|J|=n+1$. In other words we must
include all those equations which have a full row of zeros in the right
hand block. This clearly includes the equations corresponding to rows in
the bottom part of the Jacobi matrix. However, if the upper right
submatrix is not of maximal rank then we may be able to form
integrability conditions. If for a row with all zeros in the right hand
section we find that the left hand part is independent of the rows in
the lower part of the matrix then there is indeed an integrability
condition, which can be determined by performing the same operations on
the full equations $D_i \Delta_\nu = 0$.

We will call the system of equations defined by the upper right block of
the Jacobi matrix the symbol of $\mathcal{I}$, and denote this
$\mathsf{Sym}\,\mathcal{I}$:
\[
\mathsf{Sym}\,\mathcal{I} : \left\{
  \sum_{\alpha,|J|=n} \left( \rcdbpd{\Delta_\nu}{u^\alpha_J} \right) v^\alpha_J = 0
\right.
\]
where the $v^\alpha_J$ are a new set of variables, which we order in the
same way as the the $u^\alpha_J$ when displaying the symbol as a matrix.
Notice that the entries in the matrix of $\mathsf{Sym}\,\mathcal{I}$ are
the coefficients of the highest order jet variables in the equations
defining $\mathcal{I}^{(1)}$ as can be seen by comparison with the
formal derivative.

Comparison of the ranks of $\mathcal{I}$, $\mathcal{I}^{(1)}$ and
$\mathsf{Sym}\,\mathcal{I}^{(1)}$ will enable us to determine if an
integrability condition will occur on a single prolongation and
projection. There is an integrability condition if $\mathsf{rank}\,
\mathcal{I}^{(1)}_1 > \mathsf{rank}\, \mathcal{I}$, or equivalently if
$\mathsf{dim}\, \mathcal{I}^{(1)}_1 < \mathsf{dim}\, \mathcal{I}$.
Furthermore from inspection of the Jacobi matrix of $\mathcal{I}^{(1)}$
\[
\mathsf{rank}\, \mathcal{I}^{(1)}_1 = \mathsf{rank}\, \mathcal{I}^{(1)}
 -
\mathsf{rank}\, \mathsf{Sym}\,\mathcal{I}^{(1)}
\]
We can thus systematically determine if integrability conditions arise
from a single prolongation, and if necessary find the new equations.

Henceforth we will always consider the row echelon form of the symbol.
We call $\mathsf{Sym}\,\mathcal{I}$ involutive if
\[
\mathsf{rank}\,\mathsf{Sym}\,\mathcal{I}^{(1)} = \sum_k k \beta_k
\]
where $\beta_k$ is the number of rows of class $k$ in
$\mathsf{Sym}\,\mathcal{I}$.

For a row of class $k$ we call the variables $x^1,x^2,\dots,x^k$
multiplicative variables. We now consider prolonging each equation by
its multiplicative variables only. The equations obtained in this manner
will be independent as they will have distinct pivots in
$\mathsf{Sym}\,\mathcal{I}^{(n)}$. As there are $\beta_k$ equations of
class $k$ and each has $k$ multiplicative variables then this means
there will be at least $\sum k \beta_k$ independent equations of order
$n+1$ in $\mathcal{I}^{(n+1)}$. If $\mathsf{Sym}\,\mathcal{I}$ is
involutive then we obtain all the independent equations of order $n+1$
in this manner. The equations obtained from the other prolongations
required to prolong $\mathcal{I}$ to $\mathcal{I}^{(1)}$ will thus be
dependent, of lower order, or both.

The importance of involutive symbols arises from the following theorem
which provides a criterion for involution that can be tested in a
finite number of operations:

\begin{theorem}
$\mathcal{I}$ is involutive if and only if $\mathsf{Sym}\,\mathcal{I}$
is involutive and $\mathcal{I}^{(1)}_1 = \mathcal{I}$.

\end{theorem}

\section{Cartan-Kuranishi Completion}

The central theorem of the theory of involutive sytems of differential
equations is the Cartan-Kuranishi theorem:
\begin{theorem}
For every differential equation $\mathcal{I}$ there are two integers
$j$, $k$ such that $\mathcal{I}^{(k)}_j$ is an involutive equation with
the same solution space.
\end{theorem}

The Cartan-Kuranishi completion algorithm is a straightforward
application of the two previous theorems:
\begin{tabbing}
\textbf{input} $\mathcal{I}$ \\
\textbf{repeat} \\
\qquad \=\textbf{while} $\mathsf{Sym}\,\mathcal{I}$ is not involutive 
\textbf{repeat} $\mathcal{I} := \mathcal{I}^{(1)}$ \\
\>\textbf{while} $\mathcal{I} \ne \mathcal{I}^{(1)}_1$ \textbf{repeat} 
$\mathcal{I} := \mathcal{I}^{(1)}_1$ \\
\textbf{until} $\mathsf{Sym}\,\mathcal{I}$ is involutive and $\mathcal{I}
 = \mathcal{I}^{(1)}_1$ \\
\textbf{output} involutive $\mathcal{I}$
\end{tabbing}
Therefore, given a differential equation we may first complete it to an
involutive system (if it is not already involutive) and then construct a
power series solution order by order using the algorithm
described earlier.

\section{Conclusion}

Although the algorithms described in this paper may be used to construct
formal solutions to systems of partial differential equations,
they suffer from several shortcomings in practice. Firstly, there is the
problem of setting the values of the parametric derivatives. Typically,
these must be calculated from the values of functions on submanifolds of
the base space or from values of those functions at the points of a
lattice within the base space. Secondly, many terms of the power series
must be calculated to provide solutions of comparable accuracy to those
produced by discretisation schemes, and the symbolic manipulations
involved rapidly become computationally intensive as the order
increases. Thirdly, there is the problem of the convergence of the
series. A promising approach to the circumvention of these difficulties
is the use of a hybrid method that first uses a discretisation scheme to
calculate values at lattice points, uses these values to determine
approximate values of the parametric derivatives and then constructs
power series about each of the points and smoothly joins them together
using a functional interpolation scheme. This method is currently being
implemented.

Completion to involution and the construction of power series solutions
are far from the only applications of the jet bundle formalism. As
mentioned in the introduction, jet bundles also provide the natural
setting for the analysis of the symmetry groups of systems of PDEs and
of the variational symmetries of Lagrangian systems (which are linked
by Noether's theorem to conservation laws). Symmetry
analysis is also closely related to the construction of
solutions possessing specified symmetries. Unfortunately, these
fascinating and important subjects are beyond the scope of the current
paper.

\end{document}